\documentclass[11pt,a4paper]{article}
\pdfoutput=1
\usepackage{amssymb}  
\usepackage{amsthm}
\usepackage{amsmath}
\usepackage{amsfonts}
\usepackage{moreverb}
\usepackage{graphicx}
\usepackage{float}
\usepackage{enumitem}
\usepackage{caption}
\usepackage{subcaption}
\usepackage{authblk}
\title{Designing a Nonlinear Model Predictive Controller for Fault Tolerant Flight Control}
\author[1]{Rudaba Khan}
\author[2]{Paul Williams}
\author[2]{Paul Riseborough}
\author[1]{Asha Rao}
\author[1]{Robin Hill}
\affil[1]{Department of Mathematics and Geospatial Science, RMIT University, Melbourne, Australia.}
\affil[2]{BAE SYSTEMS Australia, Melbourne Australia.}

\begin{document}

\maketitle

\begin{abstract}
This paper describes the design process for developing a nonlinear model predictive controller for fault tolerant flight control.  After examining and implementing a number of numerical techniques, this paper identifies pseudospectral discretisation as the most suitable for this design. Applying the controller to a 2D robot model shows that the nonlinear controller performs much better than the linear controller, especially in the closed loop scenario. Assuming fault detection information, applying the technique to the longitudinal motion of a  generic aircraft model shows the design to be eminently suitable for flight control.
\end{abstract}

\section{Introduction}
Most research on fault tolerant control sits within the context of large manned aircraft. With regards to unmanned aerial vehicles (UAVs), the majority of the literature describes the application of fault tolerant control (FTC) to rotorcraft rather than fixed wing aircraft.  In this paper we develop a nonlinear model predictive control (NMPC) based controller for a fixed wing aircraft model suitable for the purposes of fault tolerant flight control.  The controller is first tested and analysed on a lower order 2D robot model.  Then assuming that fault information is available we successfully demonstrate for the first time the use of an NMPC based fault tolerant flight control system for a fixed wing aircraft.\\

The main characteristic of a fault tolerant control (FTC) system is its ability to automatically cope with system faults before they turn into a serious system failure.  The integration of an FTC scheme significantly increases the ability of the system to maintain overall stability in the presence of a fault \cite{zhang2008bibliographical}.  Model predictive control, a highly promising approach to re-configurable and fault tolerant control (\cite{Kale2004}, \cite{Boskovic2001}, \cite{Gopinathan1998} and \cite{maciejowski1999fault}), focuses on \textit{what to control}, instead of \textit{how to control}.  This is a subtle, but illuminating, difference allowing for the exploitation of inherent system characteristics such as nonlinearities and cross-coupling effects by the controller, rather than their influence being minimised.  This is possible due to the underlying daisy chaining capability of MPC \cite{maciejowski1998implicit}.  For example \cite{pachter2003fault}, the primary function of the rudder of an aircraft is to provide yaw or sideways control.  However, the rudder can also have some effect on the roll of the aircraft.  Therefore, in the event of the failure of an aileron actuator, the primary control surface for roll, it is still possible to execute a limited roll manoeuvre with the rudder.  This degree of fault-tolerance in the flight control system requires a suitable re-configurable architecture to be purposefully designed and implemented \cite{maciejowski1999modelling}, which is able to re-establish control, albeit with limited capacity, and execute the required manoeuvres.  The mission can then either be continued with the failed component or aborted, the primary objective being to avert a catastrophic failure or the loss of the aircraft, and to ensure that it is brought back to ground safely.  The aim of this paper is to design a controller for exactly this purpose.\\

Much of the research in model predictive control (MPC) is based on linear MPC where a linear model of the plant is utilised.  This paper develops an NMPC controller with the ultimate goal of application to fault tolerant flight control for UAVs.  The next section  (Section \ref{sec:NMPC}) will detail the workings of MPC, followed by a discussion of NMPC.  The implementation of NMPC requires the solution of an optimal control problem at each time step, and in Section \ref{sec:NMPC} several methods of achieving a solution are investigated.  In Section \ref{sec:Brach} a small selection of these optimal control methods, namely direct single shooting, direct multiple shooting and two collocation methods the first based on Euler integration and the second on pseudospectral discretisation are applied to the Brachistochrone problem, with the results aiding in the selection of the best technique to integrate into the NMPC controller design.  The results of this section showed the pseudospectral collocation method to be the best choice.  Section \ref{sec:2DRobotModel} applies a pseudospectral based linear and nonlinear MPC controller to a 2D robot model in both open loop and closed loop scenarios.  The results show that the nonlinear controller outperforms the linear controller particularly as perturbations increase and linear assumptions are violated as is the case when faults occur.  Section \ref{sec:FlightControl} applies the final NMPC controller design to the longitudinal motion of a generic aircraft model assuming fault detection information is available.  The results successfully demonstrate the capability of our controller as a fault tolerant flight control system.  Many fault tolerant control system designs require two controllers one for the nominal case and one for the fault case.  Our system requires the design of only one controller which can handle both nominal and fault cases.  Finally Section \ref{sec:conclusion} gives the analysis followed by the conclusion. 

\vspace{-6pt}

\section{Nonlinear Model Predictive Control}\label{sec:NMPC}

\vspace{-2pt}

Model predictive control (MPC), also referred to as receding horizon control, is an advanced control technology developed by practitioners within the industrial process industry that has had considerable impact on industrial process control. MPC, being capable of handling equipment and safety constraints \cite{maciejowskiBook}, allows systems to operate at or near constraints, yielding a more efficient and profitable operation.  Unlike many other control system designs, where the model of the plant is used only for design and analysis purposes, in MPC the model is an integral part of the control algorithm and is used to predict future behaviour of the plant in order to calculate the optimal control trajectory.\\

Linear MPC, where the internal model is linear, is a thoroughly researched area and is commonly used in practice.  The internal prediction model predicts the behaviour of the plant over a future prediction Horizon, $H_p$.  The idea is to select the best input that will produce the best predicted behaviour.  A number of coincidence points are placed over the horizon with the aim of bringing the predicted output as close as possible to the reference trajectory.  This is achieved by optimising a cost function, commonly a quadratic cost via quadratic programming in the case of linear MPC.  Only the first input of the calculated trajectory is applied to the plant and the prediction window slides along by the sampling time.  Once the window slides to the next time step and the calculated input is applied to the plant, the new plant states are fed back to the controller and the whole cycle begins again.  The length of the prediction window remains fixed but slides forward by one sampling interval at each step; a process referred to as the receding horizon strategy.  To reduce the computational burden a control horizon, $H_u$, can be defined which is smaller than $H_p$.  The control inputs are calculated only along the control horizon, beyond which point the value of $u$ remains constant.  Performance stability increases as the length of $H_u$ approaches the length of $H_p$.\\

Whilst NMPC uses the same structure as linear MPC the advantage of the former is the incorporation of a nonlinear process model for highly nonlinear systems. These nonlinear models are based on ``first principles'' and are obtained from an understanding of the physical nature of the system \cite{maciejowskiBook}.  Increased ability to handle the computing demands of solving nonlinear optimization problems has lead to a rise in interest in NMPC within the control community.\\  

The interest in NMPC, which began in the 90s, has been driven by the fact that today's processes need to be operated under tighter performance specifications, with more environmental and safety constraints that can only be met when process nonlinearities and constraints are explicitly considered in the controller design \cite{allgower2000nonlinear}. The major limitation of linear MPC is that plant behaviour is described by a linear dynamic model, making it unsuitable for both moderately as well as highly nonlinear processes with large operating regimes \cite{henson1998nonlinear}.  NMPC is more frequently used in the process industry because the time scales encountered are in the order of minutes, making real-time requirements less severe than in, for example, aerospace applications.\\

While NMPC has the potential to improve process operation, it poses theoretical and practical problems that are more challenging than those associated with its linear counterpart, mainly due to the nonlinear program that must be solved online at each sampling period \cite{rau2002model}.  However, the inherent robustness of NMPC that allows it to deal with input model uncertainties without taking them directly into account, is a definite advantage, which is what makes it the focus of this paper.  The next sub-section discusses optimal control techniques that can be used to solve NMPC problems.

\subsection{Problem Formulation - Optimal Control Problem}
For the NMPC methodology to be practically feasible the optimisation must be performed within the time constraints governed by the sampling period of the application \cite{Cannon2004229}.  Hence, when designing and implementing NMPC strategies,  consideration must be given to computational efficiency \cite{Cannon2004229}.  Globally optimal NMPC methods can provide benefits over local techniques and can be successfully used for online control \cite{Long2006635}.  NMPC methods are generally based around tailoring nonlinear programming algorithms \cite{Cannon2004229} to fit the structure of the online optimization or parametrising the predictions in terms of degrees of freedom.  This directly affects the size of the online optimisation problem and in turn the computational burden of the NMPC strategy.\\

Methods used to solve optimal control problems commonly fall under two categories: direct and indirect \cite{williams2005comparison}.  Direct methods have better convergence properties than indirect methods and can be used quickly to solve a number of practical trajectory optimisation problems, hence only direct methods are considered in this paper.\\ 	

The path-constrained trajectory optimisation problem as detailed in \cite{williams2005comparison}, is an area that has been heavily researched and forms an integral part of the design of an NMPC controller.  The main aim of optimal control is to determine the state and control pair that minimises a cost functional.  That is, if a state and control pair is represented by $\left\lbrace\mathbf{x}(t),\mathbf{u}(t)\right\rbrace$ then the aim is to minimise:
\begin{equation} \label{eqn:Ch3_Cost}
J = \mathcal{M}\left[\mathbf{x}(t)\right] + \int_{t_0}^{t_f} \!\left[ \mathcal{L}\left(\mathbf{x}(t),\mathbf{u}(t),t\right)\right] \, \mathrm{d}t,
\end{equation}
subject to the nonlinear state equations:
\begin{equation} \label{eqn:Ch3_Dynamical Constraints}
\mathbf{\dot{x}}\left(t\right)=f\left[\mathbf{x}\left(t\right),\mathbf{u}\left(t\right),t\right] ,
\end{equation}
the initial and terminal constraints
\begin{eqnarray} \label{eqn:Ch3_End Point Constraints}
\psi_0\left[\mathbf{x}\left(t_0\right)\right] = 0, \\
\psi_0\left[\mathbf{x}\left(t_f\right)\right] = 0,
\end{eqnarray}
the mixed state-control path constraints
\begin{equation} \label{eqn:Ch3_pathConstraints}
\mathbf{g}_L\leq\mathbf{g}\left[\mathbf{x}\left(t\right),\mathbf{u}\left(t\right),t\right]\leq \mathbf{g}_U, 
\end{equation}
and the box constraints 
\begin{equation}\label{eqn:Ch3_boxConstraints}
\mathbf{x}_L\leq\mathbf{x}\left(t\right)\leq\mathbf{x}_U,\quad\mathbf{u}_L\leq\mathbf{u}\left(t\right)\leq\mathbf{u}_U.
\end{equation}
Here: $\mathbf{x}\in\mathbb{R}^{u_x}$ are the state variables, $\mathbf{u}\in\mathbb{R}^{u_u}$ the control inputs, $t\in\mathbb{R}$ the time, $\mathcal{M}:\mathbb{R}^{n_x} \times \mathbb{R} \rightarrow \mathbb{R}$ the terminal non-integral cost (also known as Mayer component).  $\mathcal{L}:\mathbb{R}^{n_x} \times \mathbb{R}^{n_u} \times \mathbb{R} \rightarrow \mathbb{R}$ is integral cost (known as the Bolza component), $\psi_0 \in \mathbb{R}^{n_x} \times \mathbb{R} \rightarrow\mathbb{R}^{n_0}$ represents the initial point conditions, $\psi_f \in \mathbb{R}^{n_x} \times \mathbb{R} \rightarrow\mathbb{R}^{n_f}$ the final point conditions, $g_L \in \mathbb{R}^{n_x} \times \mathbb{R}^{n_u} \times \mathbb{R} \rightarrow\mathbb{R}^{n_g}$ the lower bounds on the path constraints, and $g_U \in \mathbb{R}^{n_x} \times \mathbb{R}^{n_u} \times \mathbb{R} \rightarrow\mathbb{R}^{n_g}$ the upper bounds on the path constraints.\\

Solving the problem defined by equations \eqref{eqn:Ch3_Cost} to \eqref{eqn:Ch3_boxConstraints} is difficult.  The direct methods detailed here solve this problem by applying a discretisation process and using standard algorithms to solve the resulting discrete optimisation problem.\\

A multitude of discretisation techniques exist for converting the continuous time problem to discrete time.  In direct methods the mathematical programming problem, equations \eqref{eqn:Ch3_Cost} to \eqref{eqn:Ch3_boxConstraints}, is solved by considering either discretised inputs, or a combination of discretised inputs and states, as decision variables.  The most common direct methods use the controls and states as optimisation parameters.\\

When deciding on a discretisation process many important factors need to be considered \cite{williams2005comparison}, such as the accuracy of the solution for a particular discretisation method given a number of optimisation variables, the computational expense of a particular discretisation method and the robustness of the discretisation method to the initial guess.  The methods chosen for implementation and analysis are direct single shooting \cite{Diehl2002577}, direct multiple shooting \cite{Diehl2002577}, direct collocation using Euler integration \cite{Diehl2005} and direct collocation using pseudospectral discretisation \cite{williams2011quadrature}.

\subsubsection{Shooting methods} - are used to solve the given problem with initial and terminal conditions.  They convert two point boundary value problems (BVPs) into initial value problems (IVPs) by guessing the value of the derivative at the initial boundary.  Every time a guess is made a ``shot'' (using DE solvers to find a solution), is fired in an attempt to hit the end boundary.  It is an iterative process where ``shots'' are made until the end boundary is reached within a desired tolerance.\\

The main difference between direct single and direct multiple shooting methods is that in the latter a BVP is converted into multiple IVPs.  The interval of computation, $\left[t_0,t_f\right]$ is divided into $M$ subintervals and an IVP is solved over each subinterval.  All of the solutions over the subintervals are pieced together to form a continuous trajectory/solution, whenever the solutions to the IVPs match at the beginning and end of each subinterval, the matching conditions.  These matching conditions introduce algebraic equations which must be satisfied along with the boundary conditions.

\subsubsection{Direct Transcription} - involves fully discretising the problem (all controls and all states) and then solving the discrete problem numerically.  The discretisation method, either integration or differentiation based, used to approximate the state equations must be combined with a method for approximating the integral in the generalised Bolza problem.  Pseudospectral methods, for example, are differentiation based methods relying on differentiating Lagrange Polynomial expansions of the approximating polynomials for the states, while Hermite-Simpson based techniques are often thought of as integration methods.  A comparison of various methods can be found in \cite{williams2005comparison}.  We investigated both an integration based method as well as a differentiation based method.\\

Derivative Based (or Pseudospectral) methods provide a better convergence rate known as spectral accuracy (\cite{greengard1991spectral},  \cite{ross2003legendre}).  The underlying idea is to represent the solution $f$ via a truncated series expansion and to use analytic differentiation of the series to obtain spatial derivatives of $f$.  The spectral differentiation matrix, $\mathcal{D}_N$, is a linear mapping of a vector of $N$ function values $\lbrace f\left(x_i\right)\rbrace$ to a vector of $N$ derivative values $\lbrace f^{\prime}\left(x_i\right)\rbrace$, and its calculation depends on the choice of the approximating series and the location of the points $\lbrace x_i \rbrace$.\\

One advantage pseudospectral methods have over finite element or finite difference methods is that the underlying polynomial space is spanned by orthogonal polynomials that are infinitely differentiable global functions (\cite{ross2003legendre}, \cite{elnagar1995pseudospectral}, \cite{fahroo2002direct}).  The choice of collocation points is crucial in pseudospectral methods \cite{ross2012review} and the Legendre-Gauss-Lobatto (LGL) points are implemented here because they provide maximum accuracy for quadrature approximations while at the same time avoiding the Runge phenomenon during interpolation \cite{williams2009hermite}.  In the case of LGL, the nodal points, which are the zeros of the derivatives of the Legendre polynomials, lie in the interval $\left[-1, 1\right]$ with the end points of the interval being included in the discretisation.\\  

Finally an NLP solver is required as it forms an integral part of NMPC.  Different discretisation methods are affected by the choice of NLP solvers in terms of the speed and robustness of the solution obtained \cite{williams2005comparison}.  For this work SNOPT \cite{gill2006user} is the solver of choice due to its popularity and because it is readily available.  SNOPT solves the quadratic programming subproblem with a quasi-Newton approximation to the Hessian, via a large-scale sparse sequential quadratic programming algorithm.

\vspace{-6pt}

\section{Brachistochrone}\label{sec:Brach}
The Brachistochrone problem is a nonlinear, nontrivial problem with an analytic solution very similar to the 2D robot problem to be addressed later.  The analytical solution allows determination of the accuracy of the implementation of each method and provides a benchmark in choosing a numerical method to continue development of an NMPC controller.  The Brachistochrone problem, simply stated, is to find the shape of a wire such that a bead sliding on the wire without friction, in uniform gravity, will reach a given horizontal displacement in minimum time.\\

The analytical solution is given by:
\begin{align} \label{eqn:Ch3_Brach_analytical}
x_b &= \frac{g}{\omega^2}\left(\omega\,t - \sin \omega\,t\right),&
y_b &= \frac{g}{\omega^2}\left(1 - \cos \omega\,t\right),
\end{align}
where $\omega=\sqrt{\frac{\pi \, g}{x_f}}$, $x_b$ and $y_b$ are the horizontal and vertical displacements of the bead in the $xy$-plane, $g$ is the gravitational force and $x_f$ is the final x-displacement.  The optimal control problem is to minimise the cost function:
\begin{equation}
\min \,  t_f,
\end{equation}
subject to the equations of motion of the bead:
\begin{eqnarray}
\dot{x}_b &=& V \sin \theta, \\
\dot{y}_b &=& V \cos \theta, \\
\dot{V} &=& g \cos \theta,
\end{eqnarray}
and the initial and terminal constraints:
\begin{eqnarray}
x_b\left(0\right) &=& 0,\\
y_b\left(0\right) &=& 0,\\
V\left(0\right) &=& 0,\\
x_b\left(t_f\right) &=& x_f.
\end{eqnarray}

Here $t_f$ is the time taken to reach $x_f$ and $V$ is the speed of the bead.\\

The number of discretisation points was varied for each method to investigate their effect and to determine the method most suitable for developing the fault tolerant controller.  The value of $x_f$ is set to $0.5m$ and the value of $g$ for this work is $1m/s^2$.  For the direct single shooing method the control points, $N_u$, were varied from 5 to 500 and for each value of $N_u$ the state points were varied from 10 to 1000.  Similarly the number of sections for the direct multiple shooting method was varied from $M =2$ to $M=30$ and the control points chosen for each section went from $N_u = 2$ to $N_u = 50$.  The coincidence points for both the collocation methods (Euler integration and Pseudospectral) varied from $N = 5$ to $N = 800$.\\

The accuracy of the different numerical methods was assessed by comparing the solutions produced by them and the analytical solution given in \eqref{eqn:Ch3_Brach_analytical}.  The comparison was performed by producing plots of the magnitude of the errors in $x_b$ and $y_b$.  The results showed that the Pseudospectral method produced the most accurate solution and increasing the number of coincidence points beyond $N = 50$ did not increase the accuracy of the solution.  Hence for this reason the Pseudospectral method with $N = 50$ was used as the nominal solution in the analysis of the CPU time.\\

The time taken to reach $x_f = 0.5$ in the solution of the analytical problem given by equation \eqref{eqn:Ch3_Brach_analytical} is $t_f = 1.2533\,\text{secs}$.  Plots of the CPU time taken to reach an optimal solution as a percentage of the nominal were produced.  The error plots show the percentage error between the optimal solution produced by each method and the nominal solution.\\

\begin{figure}[H]
\hspace{0.1in}
\includegraphics[scale=0.35]{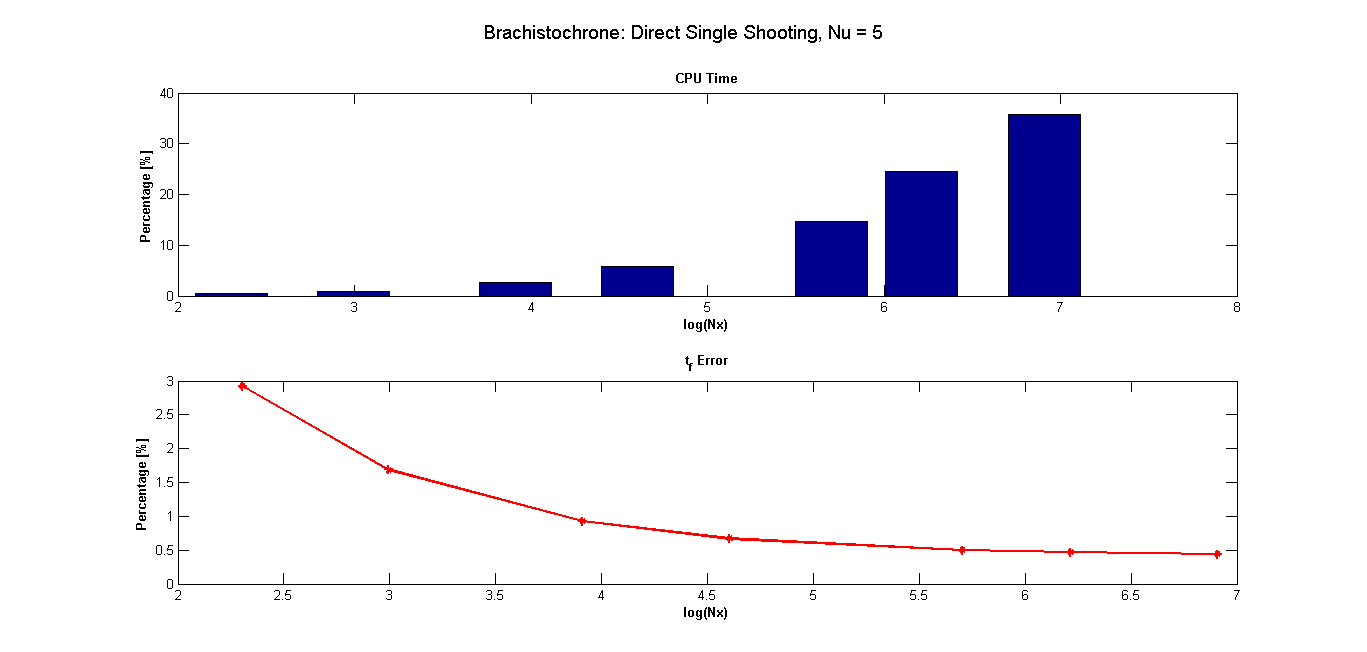}
\caption{Brachistochrone: Direct Single Shooting CPU time and $t_f$, $N_u = 5$}
\label{fig:Chap3_DSS_time_Nu5} 
\end{figure}

\begin{figure}[H]
\hspace{0.1in}
\includegraphics[scale=0.35]{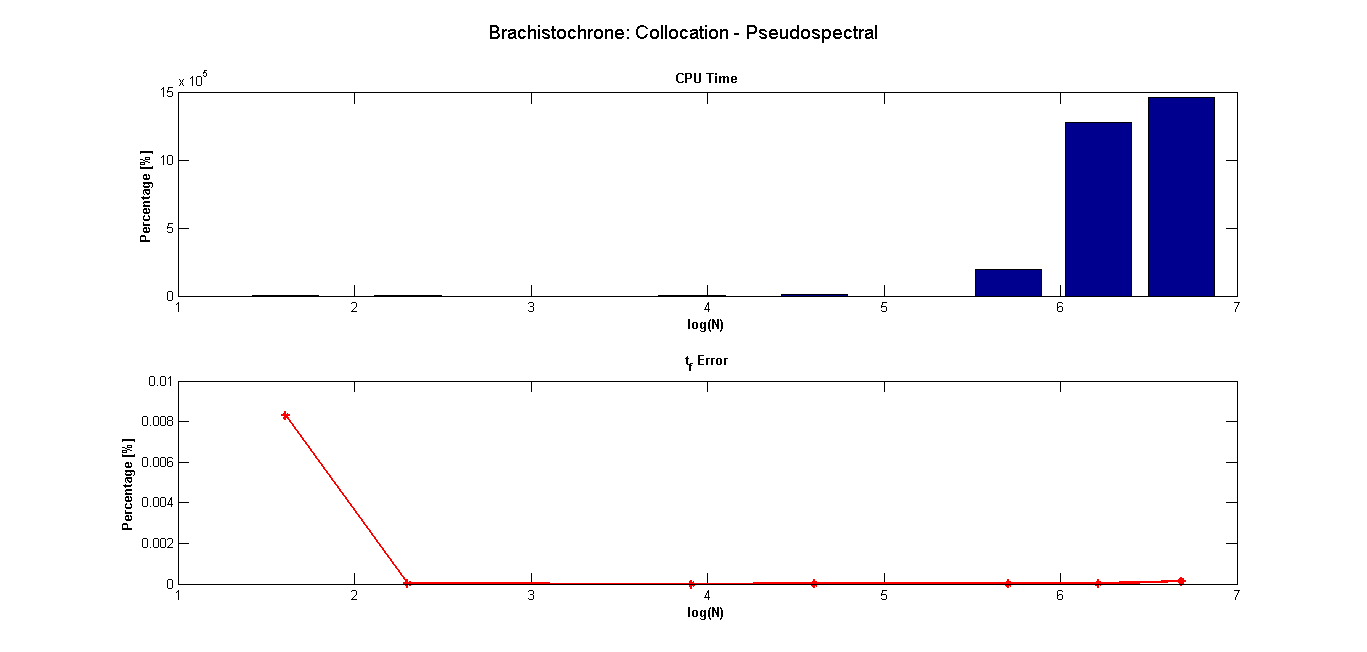}
\caption{Brachistochrone: Collocation - Pseudospectral CPU time and $t_f$}
\label{fig:Chap3_Pseud_time} 
\end{figure}

An example plot for the direct single shooting method is given in figure \ref{fig:Chap3_DSS_time_Nu5} for $N_u = 5$.  The plots showed that for $N_u = 5$ the CPU time is less than the nominal for all $N_x$ however the correct final time is unattainable.  The CPU times for the pseudospectral method are given in figure \ref{fig:Chap3_Pseud_time}.  As expected the CPU time increases for increasing $N$ and in general the CPU time taken by the pseudospectral method is higher when compared  with the other methods.\\

Overall the results showed that the pseusdospectral method can produce more accurate results with fewer discretisation points consequently requiring less time.  While, for large values of $N$, the pseudospectral method results in a greater CPU time, larger values of N are deemed unnecessary to obtain a high level of accuracy.  For this reason only the Pseudospectral method with $N=50$ points is used in design of the NMPC controller.

\section{Linear and Nonlinear MPC}\label{sec:2DRobotModel}
In this part of the design phase of the NMPC controller the pseudospectral method is integrated into an MPC framework.  An NMPC controller is designed, implemented and tested for a 2D robot model for trajectory following.  Both the open and closed loop problems are addressed and comparisons are made to linear MPC.  The next section details the 2D robot model.

\subsection{Equations of motion}\label{subsubec:chap3_Robot_EOM}
\begin{figure}[H]
\begin{center}
\includegraphics[scale=0.5]{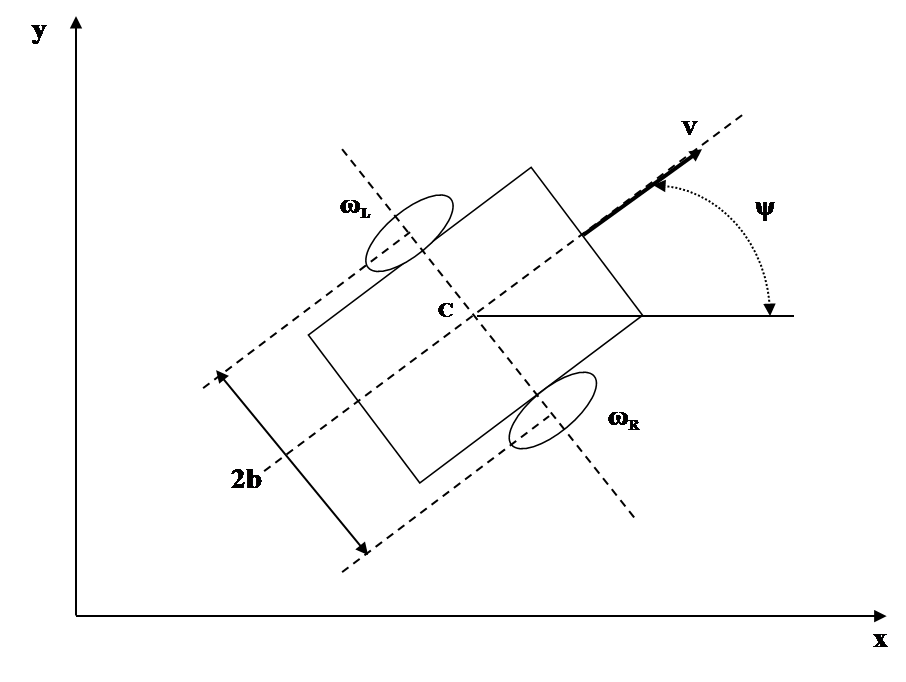} 
\caption{Robot Schematic}
\label{fig:chap3_Robot EOM}
\end{center}
\end{figure}

The 2D robot model given in figure \ref{fig:chap3_Robot EOM} is used for both the linear and nonlinear implementations of MPC.

The equations of motion for this robot model are:
\begin{eqnarray}
\dot{x} &=&  V\cos\psi, \\
\dot{y} &=&  V\sin\psi, \\
\dot{\psi} &=& \frac{R(\omega_R - \omega_L)}{2b},
\end{eqnarray}
where $x$ is the $x$-coordinate of the point $C$, $y$ the $y$-coordinate of the point $C$, $\psi$ the heading angle, $\omega_R$ the right wheel angular velocity, $\omega_L$ the left wheel angular velocity and $V$ the speed given by $V = \frac{R(\omega_R + \omega_L)}{2}$.\\

The next few sub-sections detail the development of the linear and nonlinear MPC controllers.  For a fair comparison the pseudospectral method with 50 collocation points is chosen as the method for discretisation for both controllers.

\subsection{The Open Loop Problem} \label{subsec:Chap3_OL}
In MPC an open loop problem is solved at each time step.  Hence it was important to implement and test the controller on the open loop problem.\\

Many tuning parameters can be used to determine the performance of the controller; the weighting factors on the cost function, the design of the cost function, the length of the prediction horizon, the initial condition, the integration time step, and the number of discretisation points required for an acceptable solution.  From the previous analysis the numerical technique chosen for the application of NMPC is the Pseudospectral method with 50 discretisation/coincidence points.  As a part of this work the effect of the choice of the cost function, the prediction window length, the integration time step and the effect of the initial condition on the solution were all considered.\\

The 2D robot is required to follow the path:
\begin{equation}
\forall \, x \geq 0: y = 5,
\end{equation}
travelling with a velocity of $1m/s$ and constraints of $\pm 1000 \text{deg/sec}$ on the wheel speeds $\omega_R$ and $\omega_L$.  
The objective is to drive the robot back to the reference path, from $y = 6$ to $y = 5$.

\subsubsection{Effect of Different Cost Functions}\label{para:ch3_costTypes}
Five different cost functions were developed:\\
\textbf{Cost Type 1:} Errors between the reference/nominal path and the robot path are minimised:
\begin{equation}
{J_N}_1 = \frac{\left(t_f-t_0\right)}{2}\sum_{j=0}^N\,\left(\big\Vert\mathbf{x}-\mathbf{x}_{\text{ref}}\big\Vert_{Q_x}^2\right)\,w_j.
\end{equation}
\textbf{Cost Type 2:} Errors between the robot path and the nominal path, plus the error between the actual wheel speeds, $\omega_R$ and $\omega_L$ and the nominal wheel speeds are minimised:
\begin{equation}
{J_N}_2 = \frac{\left(t_f-t_0\right)}{2}\sum_{j=0}^N\,\left(\big\Vert\mathbf{x}-\mathbf{x}_{\text{ref}}\big\Vert_{Q_x}^2 + \big\Vert\mathbf{u}-\mathbf{u}_{\text{ref}}\big\Vert_{Q_u}^2\right)\,w_j.
\end{equation}
\textbf{Cost Type 3:} Errors between the robot path and the nominal path, plus the difference between the wheel speeds are minimised:
\begin{equation}
{J_N}_3 = \frac{\left(t_f-t_0\right)}{2}\sum_{j=0}^N\,\left(\big\Vert\mathbf{x}-\mathbf{x}_{\text{ref}}\big\Vert_{Q_x}^2 + \big\Vert\omega_R-\omega_L\big\Vert_{Q_\omega}^2\right)\,w_j.
\end{equation}
\textbf{Cost Type 4:} Errors between the nominal speed and robot speed as well as the errors between the nominal angular acceleration and the robot's angular acceleration are minimised:
\begin{equation}
{J_N}_4 = \frac{\left(t_f-t_0\right)}{2}\sum_{j=0}^N\,\left(\big\Vert V-V_{\text{ref}}\big\Vert_{Q_V}^2 + \big\Vert\dot{\psi}-\dot{\psi}_{\text{ref}}\big\Vert_{Q_\psi}^2\right)\,w_j.
\end{equation}
\textbf{Cost Type 5:} Errors between the nominal speed and robot speed as well as the errors between the nominal angular acceleration and the robot's angular acceleration along with the errors between the nominal path and the robot path are minimised:
\begin{equation}
\label{eqn:J5}
{J_N}_5 = \frac{\left(t_f-t_0\right)}{2}\sum_{j=0}^N\,\left(\big\Vert\mathbf{x}-\mathbf{x}_{\text{ref}}\big\Vert_{Q_x}^2+\big\Vert V-V_{\text{ref}}\big\Vert_{Q_V}^2 + \big\Vert \dot{\psi}-\dot{\psi}_{\text{ref}}\big\Vert_{Q_\psi}^2\right)\,w_j.
\end{equation}

Each cost type was tested using both linear and nonlinear MPC.  The robot initial $x,\,y$ and $\psi$ was set to $\mathbf{x}_0 = [0\,6\,0]^\intercal$ for both controllers.  The prediction window length was varied between $H_p = 1\,\text{sec}$, $H_p = 5\,\text{secs}$ and $H_p = 10\,\text{secs}$.  Through trial and error the weights were set to: 

\[Q_x = 10, \, Q_u = 1, \, Q_{\omega} = 1, \, Q_V = 1, \, Q_{\psi} = 1\]

The optimal trajectories produced by all the different cost functions for the varying window lengths were plotted for both the linear and nonlinear MPC cases.  Results showed that for a window length of 1 sec cost types 1, 2, 3 and 5 were able to drive the robot back onto the desired path by the end of the window for both the linear and nonlinear cases, however cost type 4 was unsuccessful in doing so.\\

The error plots for all the prediction window lengths were generated showing the magnitude of the error in the y-direction between the nominal path ($y=5$) and the actual robot path.  Overall the results showed that for a path following scenario, it is best to not only minimise the path errors, but to also follow a velocity profile to obtain a smoother non oscillating solution.  For this reason cost type 5 (equation \eqref{eqn:J5}) was chosen for the final controller design and is used in the remainder of this analysis.  In addition, the error plots show that the difference in errors produced by the linear and nonlinear controllers are the least for this cost function, making it the best candidate for comparison purposes.\\

The effect of the integration time step was also investigated.  The analysis given above considers only the optimal solution produced by the controller, however in an MPC framework only the first output is applied to the plant, with the plant then providing sensor information to the navigation subsystem, for example, (on an aircraft) to calculate location and orientation information.  Hence it is important to understand the effect of the integration time step in conjunction with the optimal control input.  The integration time step was varied as follows; $ dt = [0.1,\,\,\,0.01,\,\,\,0.001]$ for varying $H_p$ lengths, namely 1 sec, 5 secs and 10 secs, and the optimal trajectories were plotted.  All results showed that the integration time step has very little effect on the results.  Upon further investigation the results showed that the smallest integration time step of 0.001secs was able to give a solution closest to the optimal.  The length of the prediction horizon was seen to have the greatest effect on the integrated solution with the integrated solution getting closer to the optimal solution as the look ahead increased. Another point to note is that the integrated output produced with the nonlinear controller more closely matched the optimal solution compared to the linear integrated output.\\ 

Varying the length of the prediction window showed that it is always best to have a longer window as this produced the lowest errors, particularly in the linear controller case.  In addition, the longer window allowed the robot to reach the nominal path more quickly.  The accuracy of the integrated solution increased as the window length increased, but, unfortunately a longer window resulted in an increase in computation time.\\

From the results obtained a window length of 5 secs was chosen for the final NMPC controller design as it was a good compromise between efficiency and accuracy.  A window length of 1 second proved to be too short to produce an accurate solution particularly in the case of the integrated solution.  While a window length of 10 secs produced an integrated solution closely matching the optimal solution in the nonlinear controller case, the solution produced by a window length of 5 secs, while not as precise, still managed to develop a solution closely resembling the optimal solution.  Hence for the rest of this research a window length of 5 seconds is used along with cost type 5 and an integration time step of 0.01secs.\\  

\begin{figure}[H]
\includegraphics[scale=0.35]{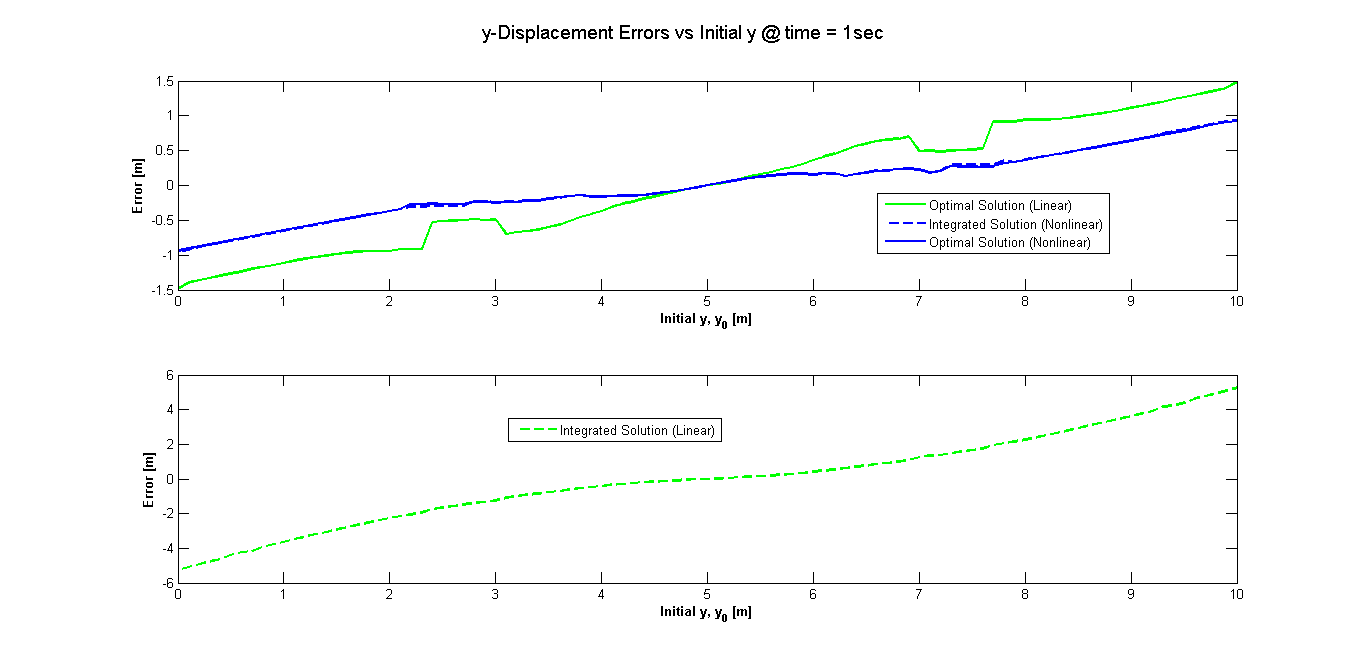} 
\centering
\caption{Open Loop: Initial Conditions vs y-Displacement Error, $\text{time} = 1\,\text{sec}$}
\label{fig:chap3_inity_1secs}
\end{figure}

\begin{figure}[H]
\centering
\includegraphics[scale=0.35]{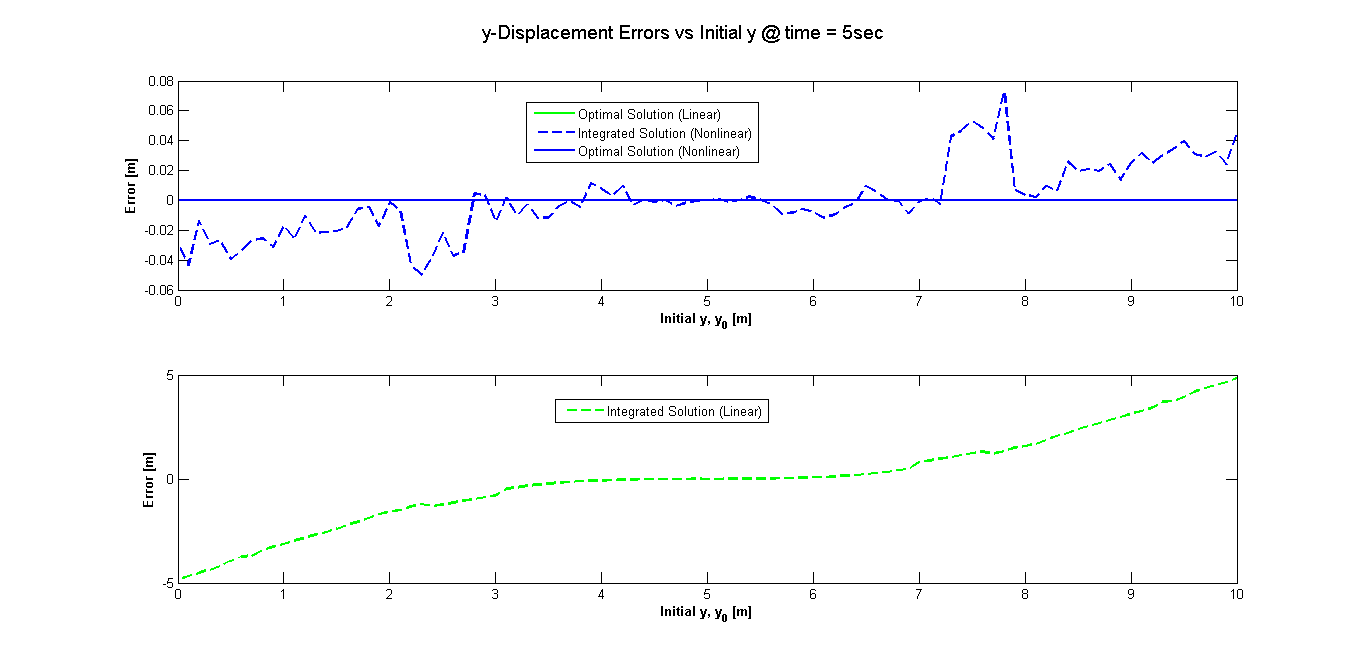} 
\caption{Open Loop: Initial Conditions vs y-Displacement Error, $\text{time} = 5\,\text{secs}$}
\label{fig:chap3_inity_5secs}
\end{figure}

The initial condition is another factor that must be considered in the design and selection of the controller.  The sensitivity of the starting point on the overall solution is critical particularly in the case of linear techniques.  In this analysis the robot is required to follow the same path as above, $y = 5$, and the initial $y$ is varied from 0m to 10m in steps of 0.1m.  The errors between the nominal path and the actual robot position are calculated at various points along the prediction horizon namely at 1 sec, 2 secs, 3 secs, 4 secs and 5 secs for all initial $y$ values.  Plots of errors versus initial $y$ for the different times were obtained and figures \ref{fig:chap3_inity_1secs} and \ref{fig:chap3_inity_5secs} show the results at $t = 1$ sec and $t = 5$ secs.  The results show the errors between the optimal solution and the nominal path as well as the errors between the integrated output and the nominal path for both the linear and nonlinear controllers.  The errors arising from the integrated output of the linear controller are shown on a separate plot underneath the main plots as these errors were much higher compared to the others and by plotting all errors on the one graph the errors produced by the other solutions were not as clearly visible.  The results show that as the time increases from 1 second to 5 secs the errors decrease as the robot approaches the nominal path.  The results clearly show that the further away the robot is from the nominal path (i.e. the greater the perturbation) the higher the error in the case of the linear controller.\\  

At the 1 second mark along the prediction window (figure \ref{fig:chap3_inity_1secs}) the errors between the solution produced by the nonlinear controller and the nominal path ($y = 5$) were seen to be linear as a function of initial $y$.  Moving further along, the prediction window shows that these errors decrease and are very close to zero for any $y_0$.  There is only a small region around the nominal path, $y = 5$, during which the errors produced by the linear controller are zero and match those produced by the nonlinear solution at any time along the prediction window.

The next subsection investigates the closed loop problem and compares the output produced by both linear and nonlinear MPC. 

\subsection{The Closed Loop Problem}
The aim of these simulations is to implement and investigate the behaviour of both linear and nonlinear MPC regarding the closed loop problem.  The fault tolerant problem is essentially closed loop, hence to apply NMPC to fault tolerant control the pseudosepectral NMPC controller design is tested on the 2D robot model where the robot is required to travel on a circular trajectory.\\

Based on the analysis from the previous subsections a prediction window length of 5 seconds is used with 50 collocation points, cost type 5 and an integration time step of 0.01secs.  Constraints of $\pm 1000 \text{deg/sec}$ are placed on the control inputs which are the angular velocities produced by the right and left wheels.  Three different scenarios were set up: Scenario 1 the robot begins on the path with initial conditions $y_0 = [5,0,0]^\intercal$, Scenario 2 the robot begins slightly off the path with initial conditions $y_0 = [-2,4,0]^\intercal$ and Scenario 3 the robot begins well off the path with initial conditions $y_0 = [0,20,0]^\intercal$.\\

For stability $H_u = H_p$ \cite{allgower2000nonlinear}.  In many MPC\textbackslash NMPC formulations $H_u$ is less than $H_p$.  While this greatly reduces the computational expense it does however produce a suboptimal solution \cite{rau2002model}, hence for the purposes of this research the control horizon is equal to the prediction horizon.

For scenario 1 where the robot begins on the path the trajectories produced by both linear and nonlinear controllers were the same.  The calculated optimal inputs produced by both controllers were also identical.  In the case of scenario 2 where the robot begins slightly off the path, both the linear and nonlinear controllers managed to bring the robot back onto the path.  The plots of the optimal inputs showed that initially both controllers work at the maximum constraint to drive the robot back onto the path.  Once the path is reached (i.e. perturbations are small) both controllers exhibit the same performance.\\

\begin{figure}[H]
\center
\includegraphics[scale=0.35]{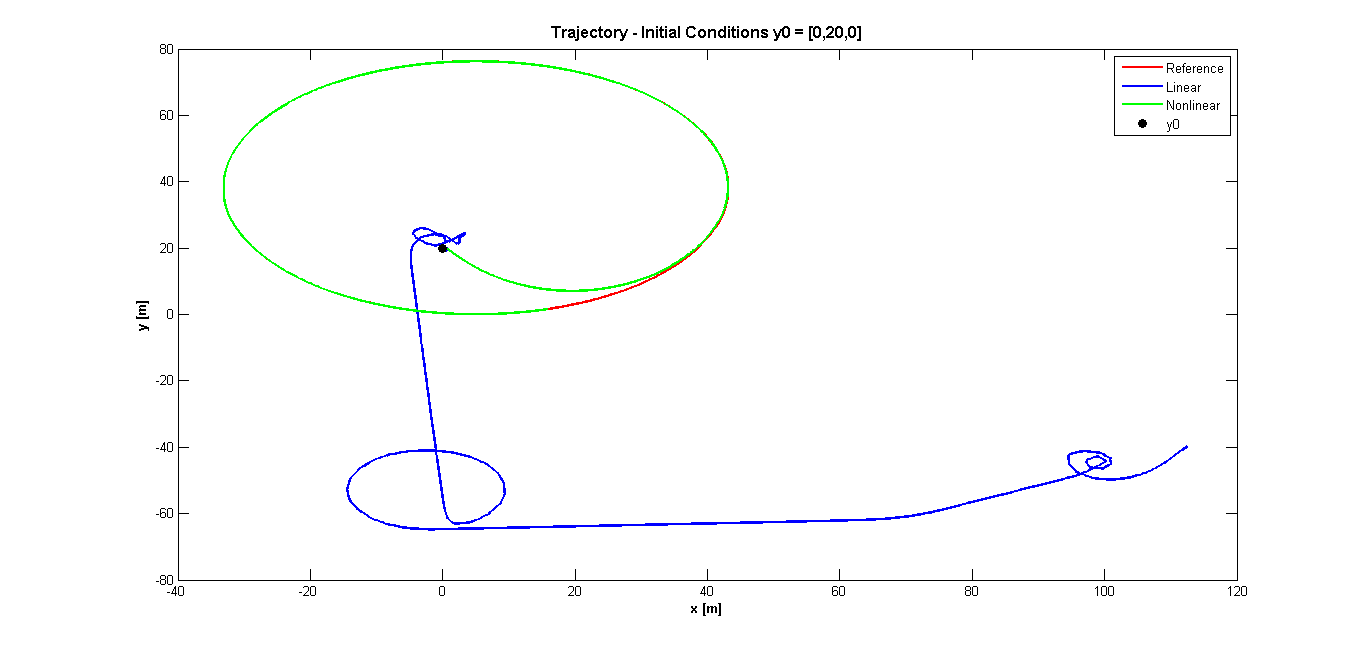} 
\caption{Closed Loop: Scenario 3 - Trajectory}
\label{fig:Chap3_CL_IC3_refTraj}
\end{figure}
\begin{figure}[H]
\center
\includegraphics[scale=0.35]{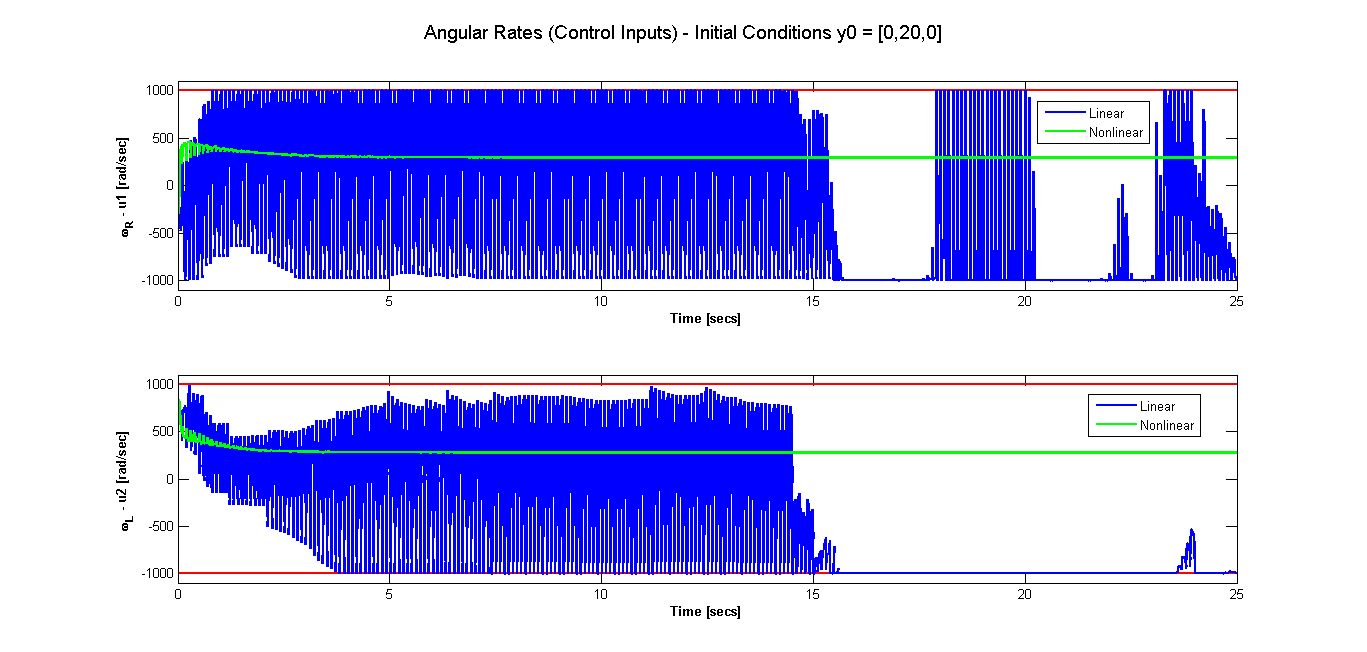} 
\caption{Closed Loop: Scenario 3 - Angular Rates}
\label{fig:Chap3_CL_IC3_omegas}
\end{figure}

The trajectory plots for scenario 3 show that only the nonlinear controller is able to bring the robot back onto the path with the linear controller unable to drive the robot back to the path.  The  control inputs produced by both controllers clearly showed that the linear controller worked very hard to take the robot back onto the path by consistently working at the constraint limits however it was still unable to return the robot back to the path.\\

The results clearly showed that the pseudospectral NMPC solution to the nonlinear model predictive controller is a viable choice outperforming its linear counterpart when the perturbations are large.  In the next section we apply the solution developed thus far to a generic aircraft model. Note that FDI is assumed for the simulation exercise carried out in section \ref{sec:FlightControl}.


\section{Application to Flight Control}\label{sec:FlightControl}
\begin{figure}[H]
\center
\includegraphics[scale=0.3]{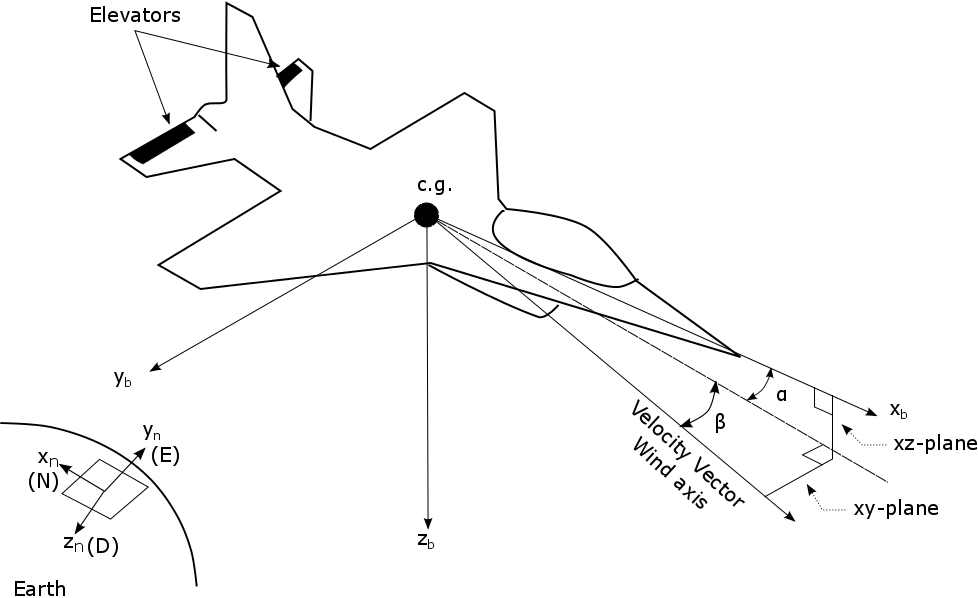}
\caption{Aircraft controls and co-ordinate systems}
\label{fig:chap5_angles_SL_paper1} 
\end{figure}

The NMPC controller developed in the previous sections was applied to the longitudinal motion of an aircraft to demonstrate fault tolerant control.  The generic aircraft model developed here for control law design and validation is based on the McDonnell Douglas F-4 aircraft \cite{garza2003collection}.  It is a fixed wing aircraft equipped with throttle, elevators, ailerons and a rudder for control.  Longitudinal motion is predominantly controlled via the throttle and elevators which is used to pitch the aircraft nose up and down and hence the remaining controls will not be considered here.  Figure \ref{fig:chap5_angles_SL_paper1} presents a sketch of a generic aircraft which identifies the location of the elevators and defines the coordinate frames in which the equations of motion are defined.  These equations are all carried out in the body axis which has its origin at the centre of gravity (c.g.) on the body of the aircraft however position and velocity information are commonly presented in the NED frame which is an earth fixed coordinate system. 
\vspace{0.1in}

The generic aircraft model has the aerodynamic characteristics of the McDonnell Douglas F-4 aircraft however the dimensional and mass properties are those given in tables \ref{table:chap5_ACDims} and \ref{table:chap5_Mass_Props_SI} respectively.   

\begin{center}
\begin{table}[H]
\caption{Aircraft Dimensional Properties}
\begin{center}
\begin{tabular}{|c|c|}
\hline 
Wing Area $S$ & 20m \\ 
\hline 
Mean Aerodynamic Chord $\bar{c}$ & 3m \\ 
\hline 
C.G location $x_{c.g}$ & 0 \\ 
\hline 
C.G reference location $x_{c.g.ref}$ & 0 \\ 
\hline 
\end{tabular}
\end{center}
\label{table:chap5_ACDims}
\end{table} 
\end{center}

\begin{table}[H]
\begin{center}
\caption{Mass Properties of model used for simulation.}
\begin{tabular}{|c|c|c|c|c|c|}
\hline 
\textbf{Parameter} & Weight (kg) & $I_X$ ($\text{kg.m}^2$) & $I_Y$ ($\text{kg.m}^2$) & $I_Z$ ($\text{kg.m}^2$) & $I_{XZ}$ ($\text{kg.m}^2$) \\ 
\hline 
\textbf{Value} & 1,177 & 2,257 & 11,044 & 12,636 & 106 \\ 
\hline 
\end{tabular}
\label{table:chap5_Mass_Props_SI}
\end{center}
\end{table}
The process model used by the NMPC controllers and the plant model is given by following equations of motion:
\vspace{-2pt}
\begin{eqnarray}
V_t &=& \sqrt{{V_N}^2 + {V_D}^2}\\
u &=& V_N\,\cos(\theta) - V_D\,\sin(\theta)\\
w &=& V_N\,\sin(\theta) + V_D\,\cos(\theta)\\
\alpha &=& \arctan\left(\frac{w}{u}\right)\\
\bar{q} &=& \frac{1}{2}\,\rho\,{V_t}^2\\
a_x &=& \frac{\bar{q}\,S\,CX + T}{m}\\
a_z &=& {\bar{q}\,S\,CZ}{m}\\
a_N &=& a_x\,\cos(\theta) + a_z\,\sin(\theta)\\
a_D &=& g - a_x\sin(\theta) + a_z\,\cos(\theta)\\
\dot{q} &=& \bar{q}\,S\,\bar{c}\,C_m\,\left(\frac{1}{I_Y}\right)
\end{eqnarray}
\begin{align}
\begin{split}
CX ={}& -0.0434 + 2.93 \times 10^{-3}\alpha + 2.53\times 10^{-5}\beta^2-1.07\times 10^{-6}\alpha \beta^2 + 9.5 \times 10^{-4}\delta_e\\
& -8.5\times 10^{-7}\delta_e \beta^2 + \left(\frac{180q\bar{c}}{\pi 2 V_t}\right)\left(8.73\times 10^{-3} + 0.001\alpha - 1.75 \times 10^{-4} \alpha^2 \right),
\end{split} \label{eqn:chap5_CX}\\
\nonumber\\
\begin{split}
C_m = {}& -6.61\times 10^{-3} - 2.67 \times 10^{-3} \alpha -6.48\times 10^{-5}\beta^2\\
& -2.65\times 10^{-6}\alpha\beta^2 - 6.54\times 10^{-3}\delta_e - 8.49\times 10^{-5}\delta_e\alpha\\
& + 3.74\times 10^{-6}\delta_e\beta^2 - 3.5\times 10^{-5}{\delta_a}^2\\
&+ \left(\frac{180q\bar{c}}{\pi 2 V_t}\right)\left(-0.0473-1.57\times 10^{-3}\alpha\right)+\left(x_{c.g.ref}-x_{c.g}\right)C_Z.
\end{split}\label{eqn:chap5_Cm}
\end{align}
Where $V_T$ is the true airspeed, $V_N$, $V_D$ are the velocities in the north and down directions in the NED frame and $u$ and $w$ are the velocities in the $x$ and $z$ directions in the body axis frame.  The accelerations $a_N$ and $a_D$ are given in the NED frame in the North and Down directions respectively and $a_x$ and $a_z$ are accelerations in the body axis.  $\dot{q}$ is the pitch rate derivative, $\bar{q}$ is known as dynamic pressure and $\alpha$ is an aerodynamic angle called the angle of attack.  $CX$ is a non-dimensional force coefficient in the body X-direction and the $C_M$ is the non-dimensional pitching moment coefficient.  The force and moment coefficients used for this model are valid for angle of attack $\alpha \leq 15 \deg$.  The thrust force, $T$ is modelled by \cite{bryson1999dynamic}:
\begin{align} \label{eqn:chap5_thrustModel}
\begin{split}
h_T = {}&\frac{H}{3048},
\end{split}\\
\nonumber \\
\begin{split}
T_{max} ={}&((30.21-0.668\,{h_T}-6.877\,{h_T}^2+1.951\,{h_T}^3-0.1512\,{h_T}^4)\\
    &+ \left(\frac{Vt}{v_s}\right)(-33.8+3.347\,{h_T}+18.13\,{h_T}^2-5.865\,{h_T}^3+0.4757\,{h_T}^4)\\
     &+\left(\frac{Vt}{v_s}\right)^2(100.8-77.56\,{h_T}+5.441\,{h_T}^2+2.864\,{h_T}^3-0.3355\,{h_T}^4)\\
      &+\left(\frac{Vt}{v_s}\right)^3(-78.99+101.4\,{h_T}-30.28\,{h_T}^2+3.236\,{h_T}^3-0.1089\,{h_T}^4)\\ 
    &+\left(\frac{Vt}{v_s}\right)^4(18.74-31.6\,{h_T}\\
    &+12.04\,{h_T}^2-1.785\,{h_T}^3+0.09417\,{h_T}^4))\frac{4448.22}{20},
\end{split}\\
\nonumber\\
T ={}& T_{max}\,\delta_{th}, 
\end{align}

The equations of motion are integrated forward in the plant model using a Runge-Kutta integration method with the Matlab subroutine ode45.  The controller runs at 10Hz and the equations of motion are used as constraints to the optimal control problem.  As developed in the previous sections, a pseudospectral discretisation method is used with 50 collocation points and a prediction window $H_p$ of 5 secs.  The optimal control inputs are calculated via SNOPT.  The aircraft is required to follow the trajectory given in figure \ref{fig:chap5_long_traj}:

\begin{figure}[H]
\begin{center}
\includegraphics[scale=0.3]{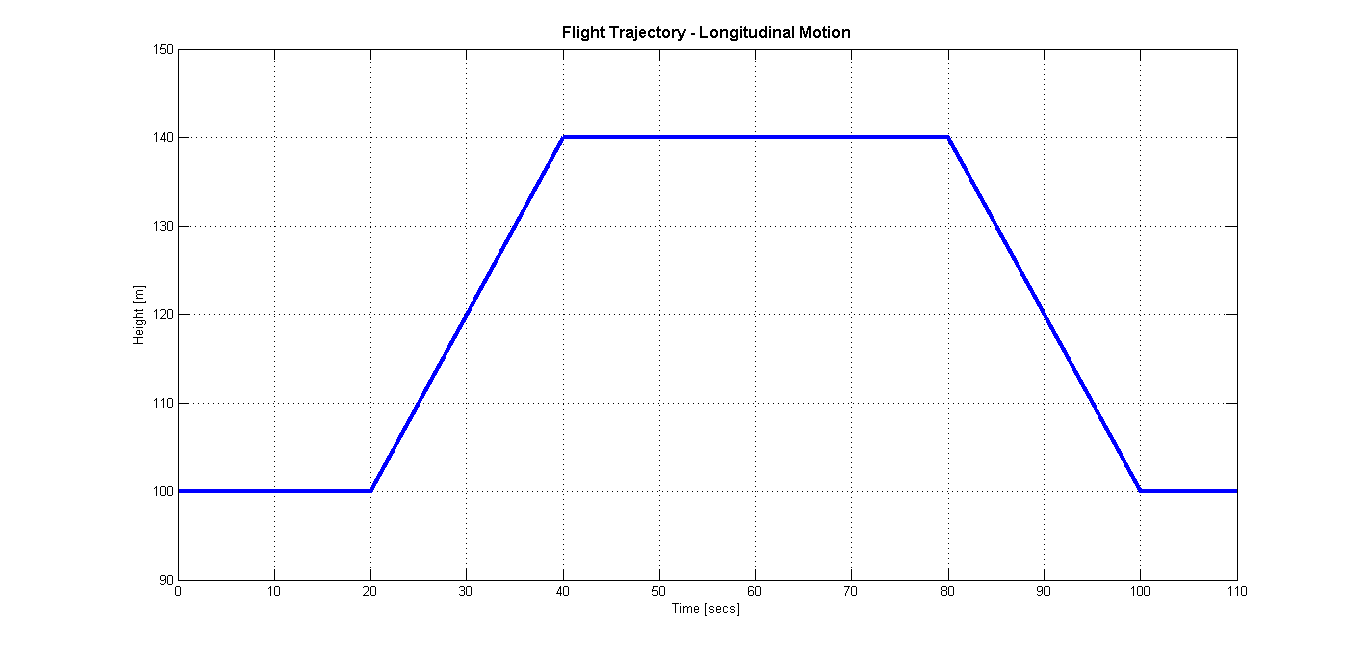}
\caption{Flight Trajectory for Longitudinal Motion}
\label{fig:chap5_long_traj} 
\end{center}
\end{figure}

Adopting the pseudospectral discretisation method where both the states and controls are discretised, the NMPC optimisation vector is:

\begin{equation}
\mathbf{x_{nmpc}} = \left[x_D,\,\,V_N,\,\,V_D,\,\,\theta,\,\,q,\,\,\delta_e,\,\,\delta_{th},\,\,\Delta{\delta}_e\right]^\intercal,
\end{equation}

where $\Delta{\delta}_e$ is the rate of change of the elevator deflection $\delta_e$.\\

The following optimal control problem is then solved:

\begin{equation}
\min_{\mathbf{x},\mathbf{u}}\, \frac{H_p}{2}\;\sum_{j = 1}^{j = N+1} \bigg(\big\Vert \bold{x}_D(j) - \bold{x}_{D_\text{ref}}(j)\big\Vert_{Q_x}^2 + \big\Vert \bold{V}_t(j) - \bold{V}_{t_\text{ref}}(j)\big\Vert_{Q_V}^2 + \big\Vert \Delta{\delta}_e\big\Vert_{Q_u}^2\bigg)\;w(j),
\end{equation}

subject to
\begin{eqnarray}
\left(\frac{t_f-t_0}{2}\right)\mathbf{D}_{j,k}\mathbf{x}_j - \mathbf{\dot{x}}_j &=& 0, \\
\mathbf{x}(j_0) - \mathbf{x}_{\text{dem}}(j_0) &=& 0,\\
\mathbf{x}_{lb}  \leq    \mathbf{x}  \leq  \: \mathbf{x}_{ub},\\
\mathbf{u}_{lb}  \leq    \mathbf{u}  \leq  \: \mathbf{u}_{ub},\\
\Delta\mathbf{\delta}_{e_{\text{lb}}}  \leq    \Delta\mathbf{\delta}_e  \leq  \: \Delta\mathbf{\delta}_{e_{\text{ub}}}, \label{eq:chap5_3DOF_cons}
\end{eqnarray}

where $\mathbf{x}_D$ and $\mathbf{x}_{D_\text{ref}}$ are the actual and reference heights respectively, and $\mathbf{V}_t$ and $\mathbf{V}_{t_\text{ref}}$ are the actual and reference true airspeeds respectively.  The constraints applied are given in table \ref{table:chap5_constraints3DoF}.

\begin{table}[H]
\begin{center}
\caption{Constraints for longitudinal Motion}
\vspace{-0.1cm}
\begin{tabular}{|c|c|c|}
\hline 
\textbf{Variable} &\textbf{ Upper Constraint} & \textbf{Lower Constraint} \\ 
\hline 
$x_D$ & 300 m & 1 m \\ 
\hline 
$V_N$ & $100\,\text{m/s}$  & $30\,\text{m/s}$ \\ 
\hline 
$V_D$ & $3\,\text{m/s}$ & $-3\,\text{m/s}$ \\ 
\hline 
$\theta$ & None & None \\ 
\hline 
$q$ & None & None \\ 
\hline 
$\delta_e$ & $20\deg$ & $-20\deg$ \\ 
\hline 
$\delta_{th}$ & $100\%$ & $0\%$ \\ 
\hline
$\Delta{\delta_e}$ & $200\,\text{deg/s}$ & $-200\,\text{deg/s}$ \\ 
\hline 
\end{tabular} 
\label{table:chap5_constraints3DoF}
\end{center}
\end{table}

The weighting matrices are diagonal matrices with the following values along the diagonal for each state which were set through trial and error:
$Q_x = 10, \, Q_V = 20, \, Q_u = 1$.  Note: the Control surface rates given in table \ref{table:chap5_constraints3DoF} are realistic for a high performance unstable airframe or for a lower weight aircraft with a stable airframe; in either case a feasible fictional aircraft model has been produced for simulation purposes to demonstrate proof of concept.

\vspace{-6pt}

\subsection{Numerical Results}
To illustrate the concept of FTC, seven different scenarios were set up, with the first one being the no fault case. The next five scenarios had the throttle stuck at $70, 50, 35, 30$ and $20\%$ for the entire duration of the flight. Scenario seven simulated the throttle getting stuck at $20\%$ 80 secs into the flight.\\

It is also assumed that fault detection information is available.  This includes the time at which the throttle becomes stuck and and the position at which it is stuck.  This information is used to update the constraint values of the NMPC controller.  Providing the controller with the most accurate and up to date information enables it to make better use of the healthy actuators.  As previously mentioned the force and moment coefficients are valid for $\alpha \leq 15 \deg$.  For all scenarios $\alpha$ was checked to ensure that $15\deg$ was never exceeded.  The plots of $\alpha$ vs time (they have not been provided here due to space constraints) showed $\alpha$ to remain below $15\deg$ hence the equations of motion were never violated.  Another means of avoiding this would be to place a constraint on $\alpha$ in the NMPC controller.   

\subsubsection{True Airspeed} 
\begin{figure}[H]
\center
\includegraphics[scale=0.3]{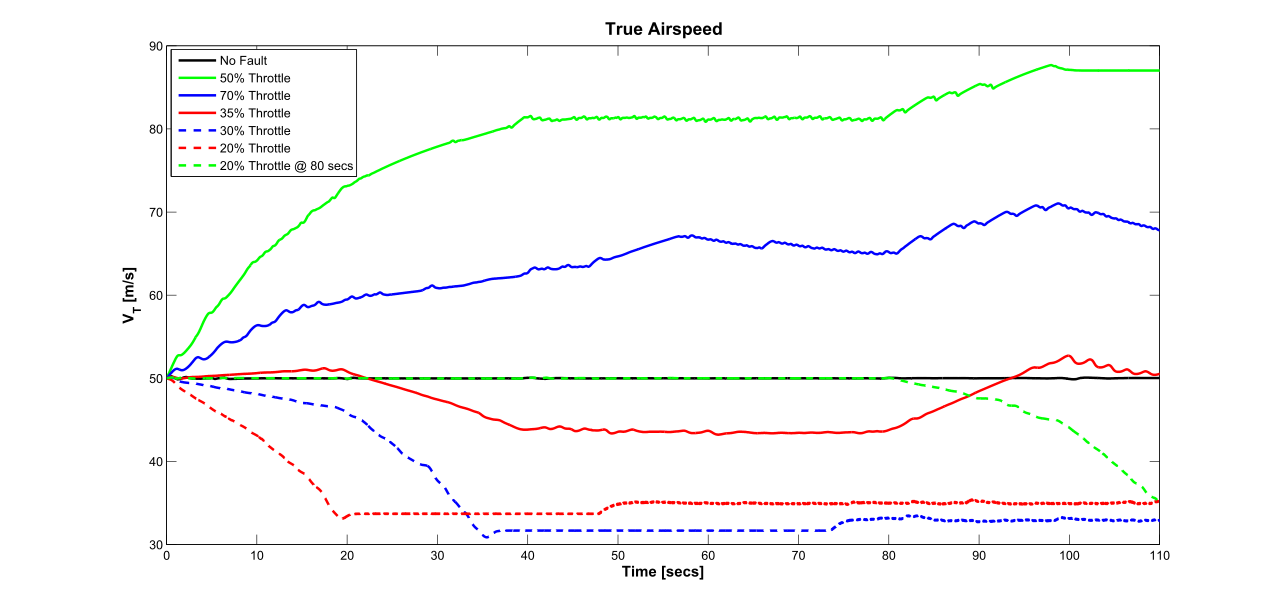}
\caption{Stuck Throttle - True Airspeed Response}
\label{fig:chap5_ST_TA}
\end{figure}

- The demanded airspeed was $50 m/s$ true airspeed.  The plot given in figure \ref{fig:chap5_ST_TA} shows the aircraft true airspeed for each of the scenarios.  The results show that in a fault free case the aircraft is able to fly at the demanded true airspeed.  However, when the throttle is stuck at $70\%$ or even $50\%$ there is too much power continually being provided to the aircraft resulting in a large airspeed response.  When the throttle is stuck at $35\%$ the aircraft is able to maintain the demanded $V_t$ for only a short period of time, at the beginning of the flight mission.  However at $30\%$ throttle the maximum deviation from the demanded airspeed is approximately $5m/s$ at any given time.  When the throttle drops below $30\%$ the aircraft is unable to maintain the true airspeed which drops to approximately between $35m/s$ and $30m/s$.  The results for scenario 7 show that once the fault occurs at 80 secs the true airspeed immediately begins to drop, as expected.  One of the main points to note is that the stall speed was never reached; the controller was able to avoid the aircraft stalling regardless of the severity of the fault.

\subsubsection{Vertical Speed} 
\begin{figure}[H]
\center
\includegraphics[scale=0.4]{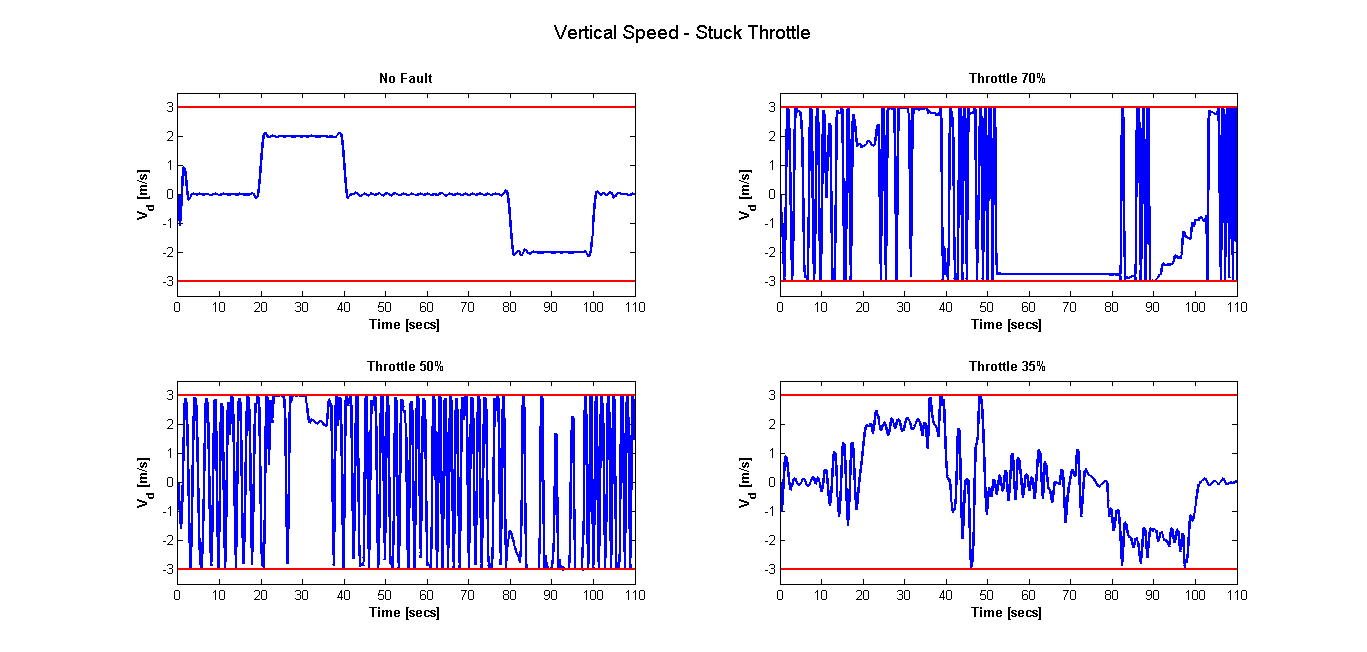}
\caption{Stuck Throttle - Vertical Speed Response, constraints (red lines), aircraft response (blue)}
\label{fig:chap5_ST_VD}
\end{figure}

- The vertical speed response of the aircraft was also analysed and plots for scenarios 1-4 are given in figure \ref{fig:chap5_ST_VD}.  The plots show the aircraft response along with the constraints (in red) placed on the vertical speed.  For high values of throttle ($70\%$ and $50\%$) the vertical speed is continuously bouncing between the constraints in an attempt to maintain the true airspeed demand.  For the case when the throttle is stuck at $35\%$ the vertical speed profile is seen to be similar to the no fault case, except in the descent phase.  During this phase, when the aircraft is descending and gaining speed, the vertical speed response can be seen to continuously move between the constraints to regulate the speed.  Results showed that for throttle values less than $30\%$ there is insufficient power to maintain a climb hence the vertical speed is seen to operate at the lower constraint or at zero.  In the case of scenario 7 it was found that once the fault occurred at 80 secs the vertical speed moved between the constraints, working hard to maintain the true airspeed.

\subsubsection{Elevator Activity} 
\begin{figure}[H]
\center
\includegraphics[scale=0.4]{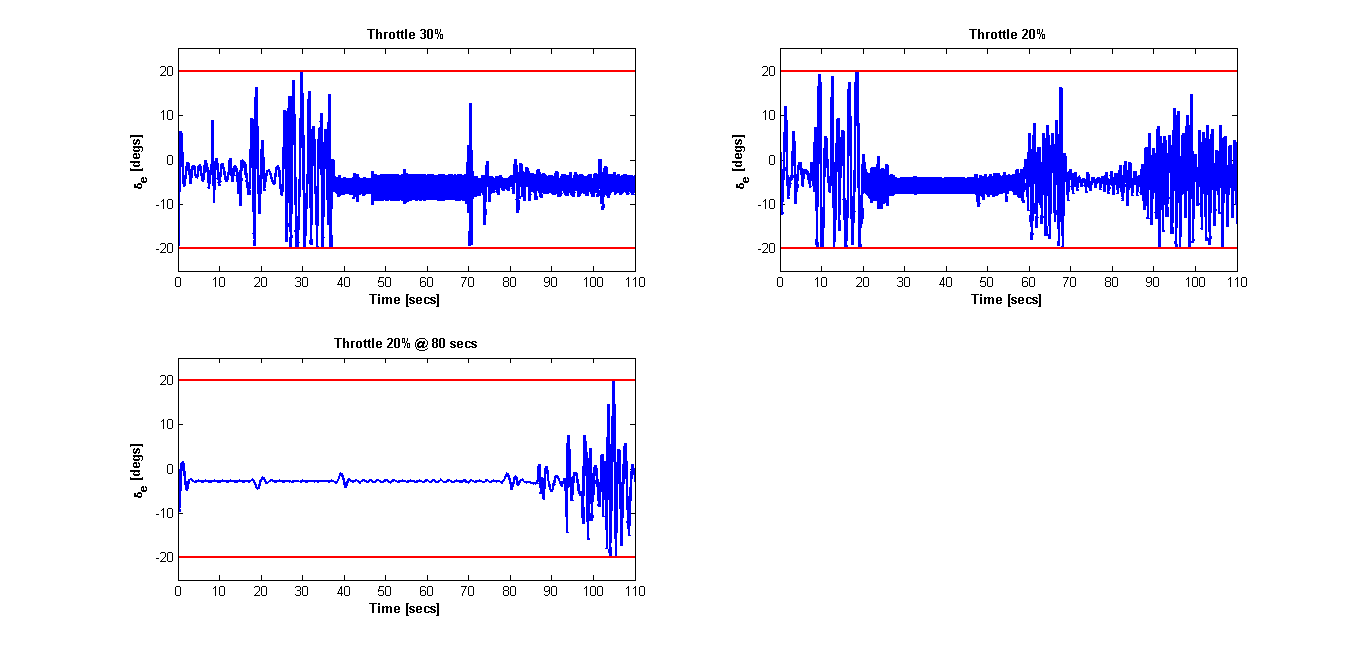}
\caption{Stuck Throttle - Elevator Activity, constraints (red lines), aircraft response (blue)}
\label{fig:chap5_ST_EA}
\end{figure}

- In regards to fault tolerance, the elevator activity is of the most interest.  If the throttle is stuck the elevator provides a level of redundancy to maintain the aircraft speed.  Figure \ref{fig:chap5_ST_EA} shows plots of elevator activity for the different scenarios.  The plots clearly show that any change in throttle increases the elevator activity when compared to the no fault case.  The elevator activity increases in an attempt to regulate the airspeed of the aircraft.  In the case of the high throttle values ($70\%$ and $50\%$) the elevator activity is the highest because a higher level of power is continually being provided to the aircraft exceeding the amount required to fly at the demanded speed.  Hence the elevator constantly jumps between the constraints in an attempt to compensate for the excess power.  For the $30\%$ stuck throttle case the elevator activity does increase compared to the no fault case; however $35\%$ throttle was found to be closer to the amount required to maintain the given height profile, hence the elevator does not need to work as hard compared to the $70\%$ and $50\%$ cases.  For the lower throttle values activity increases during the climb and descent phases.  In the climb phase of the mission there is not enough power available to the aircraft, so it compensates by erratically deflecting the elevator.  During the descent phase however there is too much power; to regulate this and to stay within the velocity constraints the activity increases.  The last scenario shows that at the fault occurrence time of 80 secs the elevator increases activity to compensate for the faulty throttle. 

\begin{figure}[H]
\center
\includegraphics[scale=0.4]{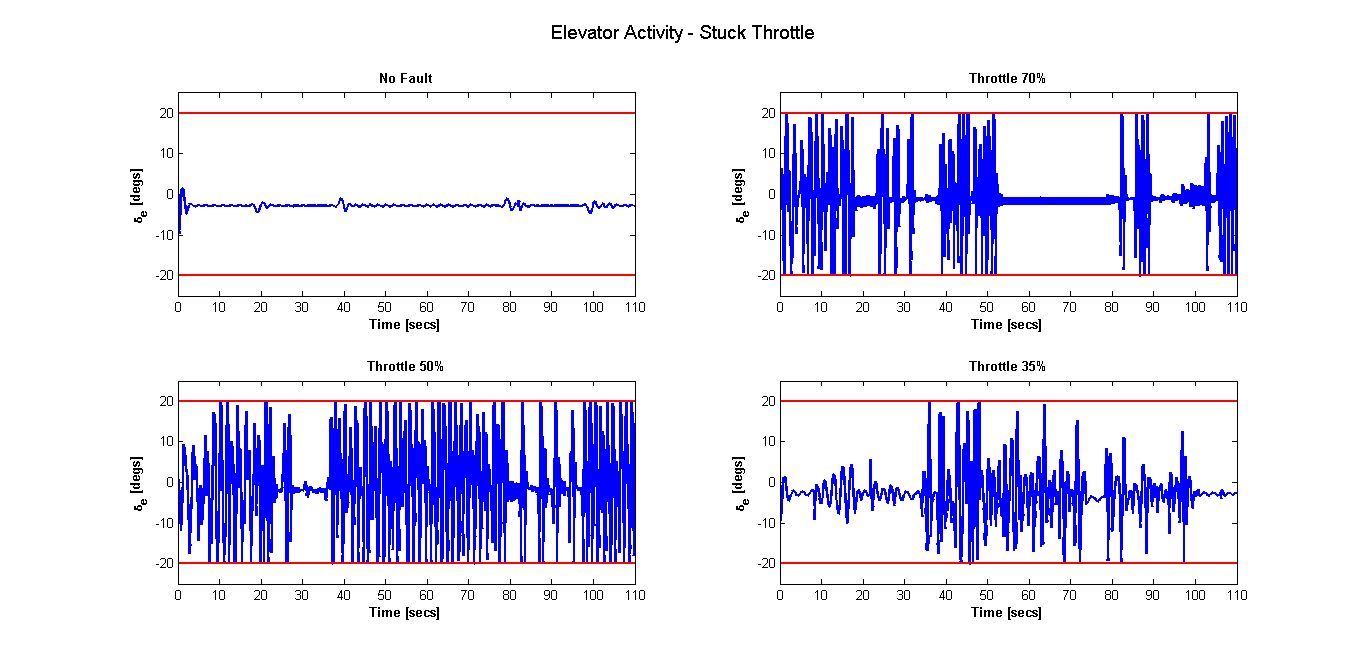}
\end{figure}
\vspace{-0.4in}

\begin{figure}[H]
\center
\includegraphics[scale=0.25]{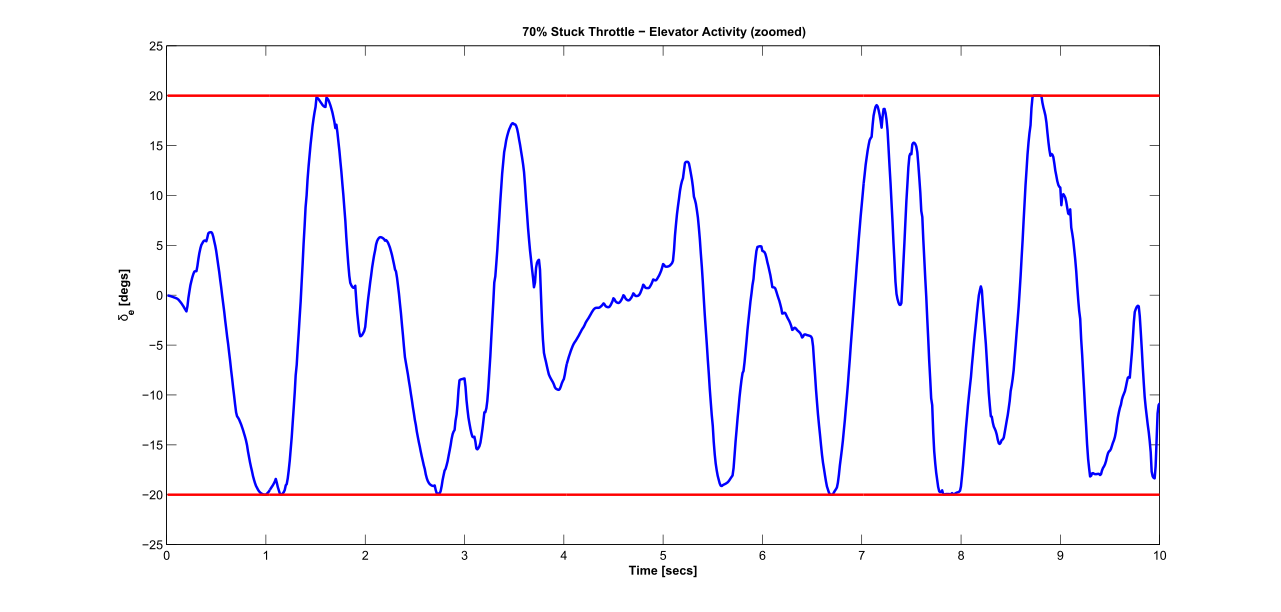}
\caption{$70\%$ Stuck Throttle - Elevator Activity (zoomed), constraints (red lines), aircraft response (blue)}
\label{fig:zoomedEA}
\end{figure}

The plots indicate large oscillations in the elevator activity however the data presented shows 200 seconds of flight.  Actuator dynamics have been modelled in both the controller and the plant model.  Zooming in on the $70\%$ stuck throttle case (figure \ref{fig:zoomedEA}) it can be seen that the rate dynamics are respected. 

\subsubsection{Height Profile} 
\begin{figure}[H]
\center
\includegraphics[scale=0.4]{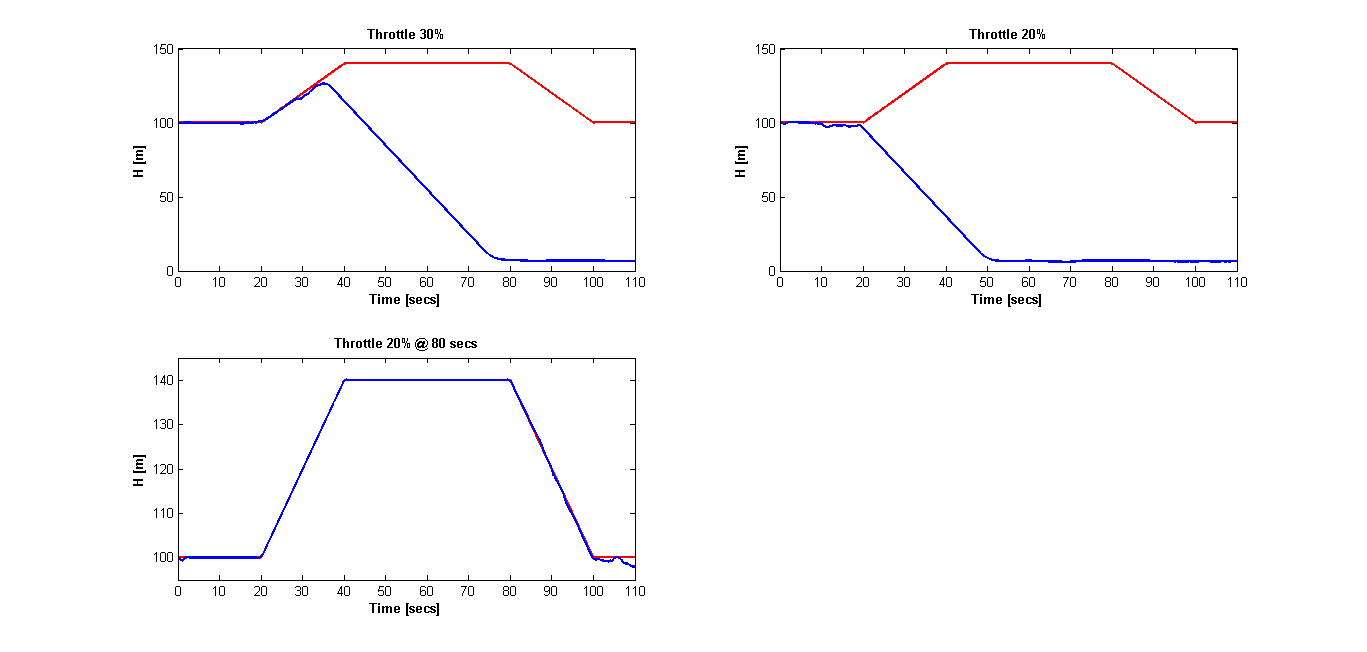}
\caption{Stuck Throttle - Height Profile, reference (red lines), aircraft response (blue)}
\label{fig:chap5_ST_HP}
\end{figure}

- The trajectory flown by the aircraft during the different scenarios is given in figure \ref{fig:chap5_ST_HP}.  The no fault case, as expected, follows the reference height profile perfectly.  The $35\%$ case is also able to closely maintain the profile.  In the $70\%$ and $50\%$ cases the aircraft continually tries to regulate the airspeed to compensate for the excess power.  The solutions produced by both scenarios show the aircraft overshooting followed by an undershoot, so the solution oscillates around the reference.  When the throttle becomes stuck at $30\%$ the aircraft begins the climb phase of the mission but is only able to continue climbing for 20 secs before it begins gliding towards the ground.  In the $20\%$ stuck throttle case the aircraft completes the straight and level phase of the mission but does not have enough power to begin climbing, and descends to the ground. The final scenario shows that the elevator is able to compensate for the stuck throttle in mid-flight and successfully finish the mission.  
\begin{figure}[H]
\center
\includegraphics[scale=0.4]{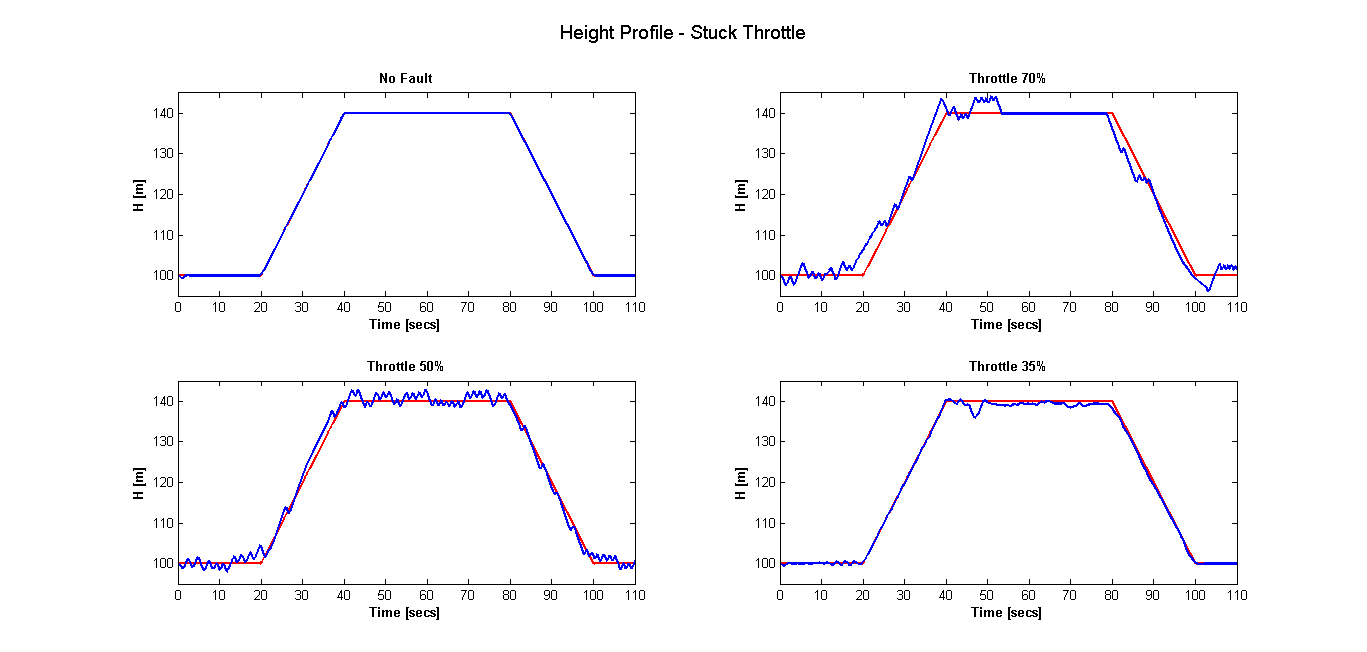}
\end{figure}
\vspace{-0.4in}

\section{Conclusion}\label{sec:conclusion} 
This paper was dedicated to the theoretical and practical implementation aspects of NMPC.  A viable design for a nonlinear MPC controller applicable to FTC was sought.  A number of discretisation methods were implemented using the well known Brachistrochrone problem and the Pseudospectral numerical method was found to have the best performance.  This method was applied in an NMPC framework to a 2D robot model in both open and closed loop settings.  Comparisons were made between linear MPC and the nonlinear MPC solutions.  For small perturbations the two controllers produced the same results.  However, for large perturbations where nonlinearity effects are more significant, the pseudospectral nonlinear MPC controller was found to produce more accurate results than the linear MPC controller.  The final design was also applied to a generic aircraft model where it was assumed that fault detection information was provided.  The motivation behind this was to demonstrate the concept of fault tolerant flight control using NMPC.  It was found that the nonlinear NMPC controller designed in this work is a viable controller for fault tolerant flight control.  This was evident from the results which show that in the event of a stuck throttle the controller is able to manipulate the movement of the elevator to compensate for the fault.


\begin{thebibliography}{9}
\bibitem{pachter2003fault} Pachter~M, Huang~Y. \emph{Fault tolerant flight control}. Journal of Guidance, Control, and Dynamics 2003; 26(1):151-160.

\bibitem{Kale2004} Kale~MM, Chipperfield~AJ. \emph{Robust and stabilized {MPC} formulations for fault tolerant and reconfigurable flight control}. Proceedings of the 2004 IEEE International Symposium on Intelligent Control 2004; 222-227.

\bibitem{Boskovic2001} Boskovic~JD, Sai-Ming~L, Mehra~RK. \emph{On-line failure detection and identification ({FDI}) and adaptive reconfigurable control ({ARC}) in aerospace applications}. Proceedings of the 2001 American Control Conference 2001; 4:2625-2626.  DOI:10.1109/ACC.2001.946269

\bibitem{Gopinathan1998} Gopinathan~M, Boskovic~JD, Mehra~RK, Rago~C. \emph{A multiple model predictive scheme for fault-tolerant flight control design}. Proceedings of the 37th IEEE Conference on Decision and Control 1998; 2:1376-1381. DOI:10.1109/CDC.1998.758477.

\bibitem{maciejowski1999fault} Maciejowski~JM. \emph{Fault-tolerant aspects of {MPC}}. IEE Two-Day Workshop on Model Predictive Control: Techniques and Applications-Day 2 (Ref. No. 1999/096) 1999; 1/1--1/4

\bibitem{maciejowski1999modelling} Maciejowski~JM. \emph{Modelling and predictive control: Enabling technologies for reconfiguration}. Annual Reviews in Control 1999; 23:13-23.

\bibitem{maciejowskiBook} Maciejowski~JM. \emph{Predictive Control: With Constraints}. Pearson Education 2002. Prentice Hall.


\bibitem{allgower2000nonlinear} Allg{\"o}wer~F, Zheng~A. \emph{Nonlinear Model Predictive Control}. Progress in Systems and Control Theory 2000; Volume 26. Birkhäuser Basel.

\bibitem{rau2002model} Rau~M, Schroder~D. \emph{Model predictive control with nonlinear state space models}. 7th International Workshop on Advanced Motion Control 2002: 136--141. {IEEE}.


\bibitem{zhang2008bibliographical} Zhang~Y, Jiang~J. \emph{Bibliographical review on reconfigurable fault-tolerant control systems}. Annual Reviews in Control, Elsevier 2008; 32(2):229--252



\bibitem{maciejowski1998implicit} Maciejowski~JM. \emph{The implicit daisy-chaining property of constrained predictive control}. Applied Mathematics and Computer Science 1998; 8:695--712.

\bibitem{henson1998nonlinear} Henson~MA. \emph{Nonlinear model predictive control: current status and future directions}. Computers and Chemical Engineering 1998; 23(2):187--202

\bibitem{Cannon2004229} Cannon~M. \emph{Efficient nonlinear model predictive control algorithms}. Annual Reviews in Control 2004; 28(2):229 - 237. DOI:http://dx.doi.org/10.1016/j.arcontrol.2004.05.001.

\bibitem{Long2006635} Long~C.E, Polisetty~PK, Gatzke~EP. \emph{Nonlinear model predictive control using deterministic global optimization}. Journal of Process Control 2006; 16(6):635 - 643. DOI:http://dx.doi.org/10.1016/j.jprocont.2005.11.001.



\bibitem{Diehl2005} Diehl~M, Bock~HG, Diedam~H, Wieber~PB. \emph{Fast Direct Multiple Shooting Algorithms for Optimal Robot Control}. Fast Motions in Biomechanics and Robotics - Lecture Notes in Control and Information Sciences 2006; 340:65--93. DOI:10.1007/978-3-540-36119-0-4.






\bibitem{Diehl2002577} Diehl~M, Bock~HG, Schlöder~JP, Findeisen~R, Nagy~Z, Allgöwer~F. \emph{Real-time optimization and nonlinear model predictive control of processes governed by differential-algebraic equations}. Journal of Process Control 2002; 12(4):577--585. DOI:http://dx.doi.org/10.1016/S0959-1524(01)00023-3.


\bibitem{williams2005comparison} Williams~P. \emph{A Comparison of Differentiation and Integration Based Direct Transcription Methods ({AAS} 05-128)}. Advances in the Astronautical Sciences American Astronautical Society by Univelt 2005; 120(1):389--409.

\bibitem{greengard1991spectral} Greengard~L. \emph{Spectral integration and two-point boundary value problems}. SIAM Journal on Numerical Analysis 1991; 28(4):1071--1080.

\bibitem{ross2003legendre} Ross~IM, Fahroo~F. \emph{Legendre pseudospectral approximations of optimal control problems}. New Trends in Nonlinear Dynamics and Control and their Applications, Springer 2003; 327--342.

\bibitem{elnagar1995pseudospectral} Elnagar~G, Kazemi~MA, R~M. \emph{The pseudospectral Legendre method for discretizing optimal control problems}. IEEE Transactions on Automatic Control 1995; 40(10):1793--1796.

\bibitem{fahroo2002direct} Fahroo~F, Ross~MI. \emph{Direct trajectory optimization by a Chebyshev pseudospectral method}. Journal of Guidance, Control, and Dynamics AIAA 2002; 25(1):160--166.

\bibitem{ross2012review} Ross~MI,  Karpenko~M. \emph{A review of pseudospectral optimal control: From theory to flight}. Annual Reviews in Control, Elsevier 2012; 36(2):182--197.

\bibitem{williams2009hermite} Williams~P. \emph{Hermite-Legendre-Gauss-Lobatto direct transcription in trajectory optimization}. Journal of Guidance, Control, and Dynamics, AIAA 2009; 32(4):1392--1395.


\bibitem{gill2006user} Gill~PE, Murray~W, Saunders~MA. \emph{Users guide for {SNOPT} version 7: Software for large-scale nonlinear programming}.   Citeseer 2006; http://citeseerx.ist.psu.edu/viewdoc/summary?doi=10.1.1.217.5929 [Last Accessed: 20-01-2016].

\bibitem{williams2011quadrature} Williams~Paul, \emph{Quadrature discretization method in tethered satellite control}.  Applied Mathematics and Computation, Elsevier  2011; 217(21):8223--8235.

\bibitem{garza2003collection} Garza~F, Morelli~EA. \emph{A collection of nonlinear aircraft simulations in {M}atlab}.  {NASA} Langley Research Center, Technical Report {NASA/TM} 2003; 212145.

\bibitem{bryson1999dynamic} Bryson~AE. \emph{Dynamic Optimization}. Addison Wesley Longman, Menlo Park, {CA} 1999.

\end{thebibliography}
\end{document}